# Compressive Sensing by White Random Convolution

Yin Xiang, Lianlin Li, Fang Li

*Abstract*—A different compressive sensing framework, convolution with white noise waveform followed by subsampling at fixed (not randomly selected) locations, is studied in this paper. We show that its recoverability for sparse signals depends on the coherence (denoted by $\mu$) between the signal representation and the Fourier basis. In particular, an *n*-dimensional signal which is *S*-sparse in such a basis can be recovered with a probability exceeding $1-\delta$ from any fixed $m \sim O(\mu^2 S \log(n/\delta)^{3/2})$ output samples of the random convolution.

## I. INTRODUCTION

THE SPARSITY characteristic of many interesting natural signals continues to be a hot topic in signal processing because it is useful for reducing signal storage requirements. Researchers in the emerging field of compressive sensing have made a step forward by relating a signal's sparsity structure with its acquisition and have showed that one can recover a sparse signal from a few linear non-adaptive measurements [1], [2], [3], [15], [16]. These measurements are realized by projecting the signal onto a series of test waveforms. For instance, if a discreted signal $x$ is a n-dimensional vector and $\{\phi_1, \ldots, \phi_m\}$ is $m$ test vectors with the same length as $x$, then the measured data are given by

$$y_i = \langle x, \phi_i \rangle, i = 1, \ldots, m \quad (1.1),$$

where $\langle \ \rangle$ stands for the inner product of two vectors. If we construct a matrix $\Phi$ from these waveforms as the row vectors, then the data equations of (1.1) can be written as

$$y = \Phi x, y = [y_1, \ldots, y_m]^T \quad (1.2),$$

where $y$ is the data vector and $\Phi$ is often called the sensing matrix.

Sometimes $x$ itself is sparse, but in this paper we consider a more general case. We assume that $x$ is *S*-sparse in an orthobasis $\Psi$. That is, if we decompose $x$ in $\Psi$ as $x = \Psi \alpha$, $\alpha$ has no more than $S$ non-zero entries. We still focus on $\ell_1$-*minimization* for signal recovery, where the estimation of $\alpha$ is given by solving the following convex optimization problem

$$\hat{\alpha} = \arg\min_{\alpha} \|\alpha\|_{\ell_1} \quad s.t. \quad y = \Phi \Psi \alpha \quad (1.3).$$

The sensing matrix $\Phi$ is often chosen as a random matrix for it keeps small coherence with $\Psi$ [2], [9]. One typical sensing system of this case is the well-known single pixel camera [7]. However, it is also known that the aforementioned random structure is not efficient for both realization and computation, especially for large-scale problems. So far, many sensing matrices with certain structures have been investigated by many authors [6], [12], [25] and it is noted that the random convolution-based approach is the most popular one due to its outstanding advantages in the aspect of hardware design and algorithm consideration.

In this paper we derive CS theorems for convolution with white noise waveforms followed by deterministic subsampling. Performance of this framework is quite simple, straightforward and not costly. Only small adjustments are required to transform a traditional sensing system to a compressive sensing one. Practical examples include the famous coded aperture [20], radar system [11], [12], [21], [22], CMOS compressed imaging [8] and so on. The main contribution of our work is giving a theoretical bound estimation for the measurement number and the requirement for signal representations for the white random convolution system.

### 1.1 White random convolution

Assume the *n*-dimensional signal $x$ is first convolved with a given white random vector $\boldsymbol{h} = [h_1, \ldots, h_n]^T$ and then subsampled. We only consider circular convolution in this paper, so $\boldsymbol{h}$ is also an *n*-dimensional vector. Since $\boldsymbol{h}$ is white, the expect vector and autocorrelation matrix of $\boldsymbol{h}$ are as follows:

$$\mu_h = E(\boldsymbol{h}) = 0 \quad (1.4),$$

$$R_h = E(\boldsymbol{h}\boldsymbol{h}^*) = \sigma^2 I \quad (1.5),$$

where $*$ denotes the conjugate operator and suppose $\sigma^2 = 1$. In this paper we focus on two typical distributions, the Gaussian distribution and symmetrical Bernoulli distribution, both of which meet the above assumptions.

Denoting the convolution matrix as $H$, it consists of the random vector $\boldsymbol{h}$ and its circular shifts, and can be written as:

Yin Xiang and Fang Li are with the Institute of Electronics, Chinese Academy of Sciences, Beijing, 100190, China. (086-010-58887481; e-mail: fli@ mail.ie.ac.cn).
Lianli Li, was with the Institute of Electronics, Chinese Academy of Sciences, Beijing, China. He is now with the Department of Petroleum Engineering, Texas A&M University, USA (e-mail: lianlinli1980@gmail.com).

$$H = \begin{bmatrix} h_1 & h_n & \cdots & h_2 \\ h_2 & h_1 & & \vdots \\ \vdots & & \ddots & h_n \\ h_n & h_{n-1} & \cdots & h_1 \end{bmatrix} \quad (1.6).$$

Without loss of generalization, we rearrange the rows of $H$ as:

$$H = \begin{bmatrix} h_1 & h_2 & \cdots & h_n \\ h_n & h_1 & & \vdots \\ \vdots & & \ddots & h_2 \\ h_2 & h_3 & \cdots & h_1 \end{bmatrix} \quad (1.7).$$

Then the $k$th row of $H$ denoted by $H^k$ is given by

$$H^k = \boldsymbol{h}^* (D)^{k-1} \quad (1.8),$$

where

$$D = \begin{bmatrix} & 1 & & \\ & & \ddots & \\ & & & 1 \\ 1 & & & \end{bmatrix}_{n \times n} \quad (1.9),$$

is the row shifting matrix.

1.2 Continuous case

For practical consideration, we assume to sense a continuous signal with limited control parameters. The signal is band-limited or is approximated by a limited-resolution observation. For simplicity, we suppose the bandwidth of the signal to equal 1. Since $y(t) = h(t) * x(t)$, the full sampled data can be given by

$$y(k) = \int h(k - \tau) x(\tau) d\tau$$

For $x(t) = \sum_{k=1} x(k) \sin c(\pi(t-k))$, we are interested in recovery of the Nyquist samples $x(k)$.

Denote $\hat{h}(k) = h(t) * \sin c(\pi t)\big|_{t=k}$, then

$$y(k) = \sum_l x(l) \int h(k-\tau) \sin c(\pi(\tau-l)) d\tau$$
$$= \sum_l x(l) \int h(k-l-\tau) \sin c(\pi \tau) d\tau$$
$$= \sum_l x(l) \hat{h}(k-l)$$

It comes back to the discrete framework of convolution. However, the entries of the random sequence $\hat{h}(k), k = 1, ..., n$ are not independent of one another in general. Since

$$E(\hat{h}(k), \hat{h}(l)) = E\begin{pmatrix} \int h(\tau) \sin c(\pi(k-\tau)) d\tau \\ \cdot \int h(\tau') \sin c(\pi(l-\tau')) d\tau' \end{pmatrix}$$
$$= \iint \delta(\tau - \tau') \sin c(\pi(k-\tau)) \sin c(\pi(l-\tau')) d\tau d\tau'$$
$$= \int \sin c(\pi(k-\tau)) \sin c(\pi(l-\tau)) d\tau$$
$$= \delta_{k,l}$$

where $\delta_{k,l}$ is the Dirac function, we can still demonstrate the continuous problem by our discrete model, $\hat{h}(k), k = 1, ..., n$, which is independent when generated from a Gaussian distribution. Accordingly, we can extend our convolution model for the continuous case for a Gaussian distribution.

1.3 Main results

In this paper we show that convolution with a white noise waveform followed by subsampling at fixed locations is an efficient CS framework to capture sparse signals. As the convolution matrix is not strictly orthogonal, we used the coherence parameter $\mu$ to measure the coherence between the sparse representation $\boldsymbol{\Psi}$ and the Fourier basis $\boldsymbol{F}$, rather than the coherence between $\boldsymbol{\Psi}$ and convolution matrix $H$ [4], [6], i.e., it is still an important parameter affecting the exact recovery.

The bound for the measurement number required to ensure exact recovery is given by the following theorems.

For Gaussian ensembles,
**Theorem 1.1**: *Suppose $\boldsymbol{\Psi}$ is an orthobasis, $h$ is an $n$-directional Gaussian white noise waveform, the convolution matrix $H$ is generated by $h$ and is shifted according to the description in 1.1. Fix a support set $T$ of size $|T| = S$ in the $\Psi$ domain, and choose a sign sequence $z$ on $T$ uniformly at random. Let $x$ be the test signal which is supported on $T$ with signs $z$ in $\Psi$, and choose samples on fixed locations $\Omega$ of size $|\Omega|=m$. Suppose that*

$$m > K\mu^2 S \log(n/\delta)^{1/2}$$
$$\cdot \max\left(\log\left(2e^2(S+1)^2\right), \log(n/\delta)\right)$$

*where $F$ is the discrete Fourier matrix, $\mu(F, \Psi) = \max_{1 \le i, j \le n} \left|(F\Psi)^i_j\right|$ is the mutual coherence between $\Psi$ and $F$, and $K$ is a numerical constant. Then with probability exceeding $1-\delta$, every signal $x_0$ supported on $T$ with signs matching $z$ can be recovered from $y = U^\Omega x_0$ by solving (1.3).*

For symmetrical Bernoulli ensembles,
**Theorem 1.2**: *Suppose $\Psi$ is an orthobasis, $h$ is an $n$-directional symmetrical Bernoulli white noise waveform,*

*the convolution matrix H is generated by h and its shifts according to the description in 1.1. Fix a support set T of size |T| = S in the Ψ domain, and choose a sign sequence z on T uniformly at random. Let x be the test signal which is supported on T with signs z in Ψ, and choose samples on fixed locations Ω of size |Ω|=m. Suppose that*

$$m > KS\mu^2(F,\Psi)\log(n/\delta)^{3/2}$$
$$\cdot \max\left(\log(4e^2(S+1)^2), \log(n/\delta)\right),$$

*where F is the discrete Fourier matrix, $\mu(F,\Psi) = \max_{1 \leq i,j \leq n} |(F\Psi)^i_j|$ is the mutual coherence between Ψ and F, and K is a numerical constant. Then with probability exceeding 1−δ, every signal $x_0$ supported on T with signs matching z can be recovered from $y = U^\Omega x_0$ by solving (1.3).*

1.4 Related works

Application of a random filter for compressive sensing was first mentioned by J. Tropp et al. [10] who proposed two equivalent realization structures of a random filter: 1) convolution with a random waveform in the time domain, and 2) multiplication with random weights in the frequency domain, both followed by equal interval down-sampling. The recovery performance for the random filter was studied with different lengths by numerical simulation. In this paper we focus on deriving the theoretical bound on the number of samples for exact recovery of sparse signals of the first structure. It should also be mentioned that Tropp proposed a dual structure for the random filter named the random demodular for efficiently sensing frequency -sparse signals [26].

Compared to J. Romberg's work [6], our work shows two significant points. In [6], the randomness is designed in the frequency domain, where the spectrum of the random waveform has unit amplitude and independent random phase such that the random waveform is orthogonal with its shift which makes the convolution matrix orthogonal. Following [4], if the convolution matrix is *orthogonal* and the sensing matrix is constructed by randomly selecting rows of the convolution matrix, the theoretical bound of the measurements number can be more readily determined by

$$m \sim O(\mu^2 S \log n)$$

where μ is the coherence between the convolution matrix and signal representation basis. In contrast, in our model the convolution matrix is not orthogonal and some frequency entries of the interested signal will be filtered. Accordingly, our model is not suitable for sensing signals which are sparse in the frequency domain. However, this sacrifice leads to an advantage of system realization. We will show that the suitability of the white convolution system for sensing a sparse signal depends on the coherence between the signal representation and the Fourier basis. Another significant difference is that we show the randomly selected sampling strategy is not necessary. We will prove that subsampling at arbitrary fixed locations also works well for the random convolution framework. However, the determinist subsampling framework for random convolution has been mentioned by H. Rauhut [28].

In [28] H. Rauhut mainly focused on sensing and recovering a time domain sparse signal by using arbitrary subset of rows of a random circulant or Toeplitz matrix. He improved the estimation of *m* given in [12] as

$$m \sim O\left(S\log(n/\delta)^3\right),$$

where *m* is the necessary number of measurements ensuring exact recovery via $\ell_1$-*minimization*, and *n*, *S*, and *δ* are of the same meaning as mentioned above. If choosing the signal representation as the identity matrix (*I*), in this case the coherence $\mu = \mu(F,I) = 1$, one can derive a similar result from the main theorems given in section 1.3. However, compared to [28] the recovery property for a general signal representation Ψ is fully studied and the starting point of the proof is completely different in the present paper.

## II. APPLICATIONS

We aim at the introduction of a random convolution framework which is close to the practical convolution system, such as SAR and coded aperture, where the independence randomization is implemented in time rather than in the frequency domain. Derivations in [2], [3], [4] show that independence of randomness plays a key role in affecting the recovery property of a random projection sensing system. Different freedom of independence and different implementation of independence result in totally different recovery properties.

In this section we describe two traditional imaging systems: SAR and coded aperture, which can both be easily transformed to a CS imaging system. Our convolution framework roughly matches these applications and is more precise than Romberg's framework.

2.1 Coded Aperture

Coded aperture is a traditional imaging system for which most current research is focused on designing the code mask and properties related to the point spread function which is used to reconstruct the original image from coded observations by linear recovery methods. Recently this old imaging framework is studied within the context of CS. Coded aperture works as a spatial convolution system, where measured data are gathered from an image convolved with the coded mask. Denote I(x1, x2) as the image scene and h as the point-spread function of the coded mask such that the coded image is given by: y=h*I. Since regular images always have a low dimensional structure, for example, sparsity in the wavelet domain, the original image can be reconstructed from a low resolution observation of y. More details can be found in [20].

2.2 SAR

SAR is a widely used remote sensing system which aims to capture the reflectivity of the target scene. As shown in Fig1, the measured SAR data is usually a two dimensional block that consists of the convolution output in both range and azimuth direction [19]. Let I(x1, x2) denote the reflectivity function of the target scene, and Rr, Ra denote the integral of I(x1, x2) along range and azimuth, respectively. If the radar is far away from the target scene (which is almost true in practice), the receiving position is fixed, the received signal is the transmitting waveform convolved with Rr, the transfer frequency is fixed, and the received signal along the aperture is the free space green function convolved with Ra(t), which is fundamental in the simplified model of SAR.

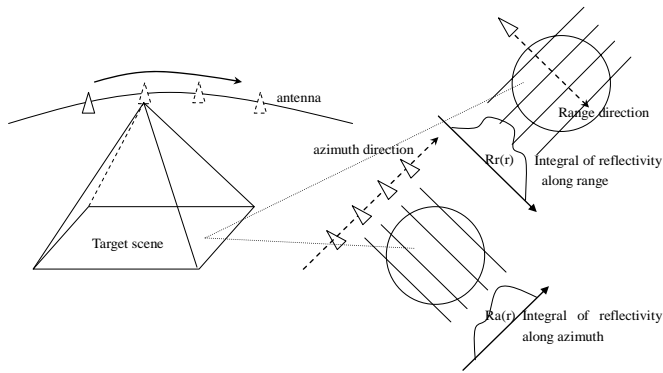

Fig. 1 SAR working strategy

The reflectivity image is sparser by taking the range and azimuth direction into account rather than one of these alone. The reduction of sampling rate in the range direction, which results in the reduction of ADC transfer speed, has been reported in [12], [17], [21], where both the random waveform and traditional frequency modulated continuous wave (FMCW) (pseudo-random waveform) are efficient for recovering the reflectivity image. For some special applications the reduction of sampling rate in the azimuth direction is more important. For example, in the application of ionosphere observation, the pulse repetition duration is not long enough to support the signal directly reflected by the ionosphere and the signal passing though the ionosphere and reflected by ground if sampled at the Nyquist rate. However, the sampling rate in the azimuth direction can be reduced for receiving both of the reflected signals, so more information is gathered from the ionosphere to engage a high resolution observation.

We simulate the SAR returns with 30-dB noise from the synthetic scene of Fig2. (a). The white area is where the object is located, the black area is its shadow, and the grey area is the uniform background. The simulation data are generated as a P-band SAR flying at 600km away from the ground, working at 435MHz with 6MHz bandwidth. The reconstruction of the target image using the conventional SAR method with the fully sampled data and using the compressive sensing recovery method with one-fourth downsampling data in both the range and azimuth directions are shown in Fig2. (b) and Fig2. (c). The CS reconstruction provides a much better result than the conventional reconstruction.

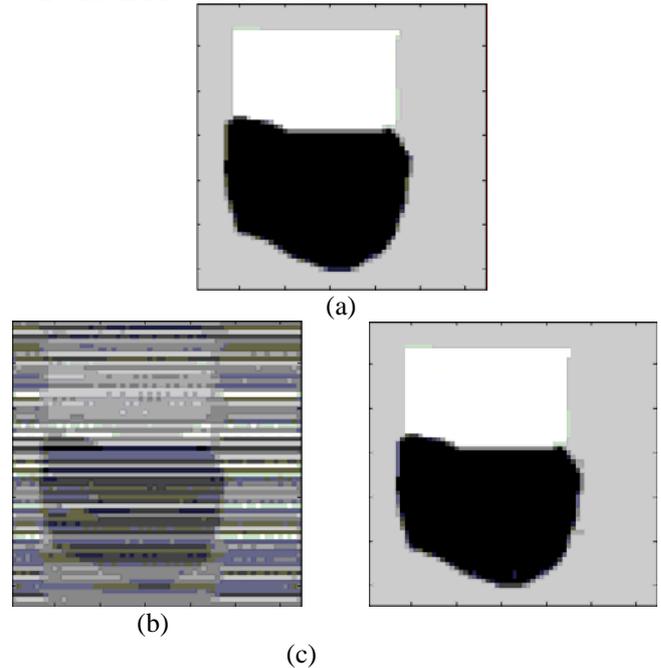

Fig2. SAR recovery from downsampled data: (a) Original scene (b) Conventional SAR method (c) compressive sensing recovery method

III. THEORY

3.1 RIP verification for white random convolution system followed by fixed down sampling

The verification of the restricted isometry property for a random convolution matrix follows the guidelines in [14], [27]. First, a new concentration inequality, in particular lemma 3.1, holds for an arbitrarily chosen (but fixed) vector with given support. The proof is a simple application of the concentration property of Lipschitz functions defined in product space, which gives an estimation of the probability tail bound in (3.2). However, this bound is dependent upon the particular choice of the vector. In the following derivation, an upper bound of the probability tail bound, in (3.1), is then developed. Secondly, the result of lemma 3.1 is generalized to any signal vector with the same support in theorem 3.2 which gives the RIP verification.

**Lemma 3.1**: *Fix an n-dim orthobasis $\Psi$, generate a random waveform $h \in R^n$ with all entries independent copies of a Gaussian random variable following the distribution of $N(0,1)$, and compose the random* *convolution matrix H with **h** as described in section 1. Let U=H**Ψ**, fix a subset Ω of size m on measurement domain, let $H^\Omega$ be the submatrix obtained by retaining the rows of H indexed by Ω, and fix a subset T of size S on signal domain. Then for an arbitrary but fixed signal $x_0$ with support on T,*

$$\Pr\left( \left| \left\| \frac{1}{\sqrt{m}} U^\Omega x_0 \right\|^2 - \|x_0\|^2 \right| > r \|x_0\|^2 \right)$$

$$< e^2 \exp\left( -\frac{mr}{2\mu^2(F,\Psi)S} \right), r > \frac{2\mu^2(F,\Psi)S}{m}$$

(3.1)

*where F is the discrete Fourier matrix, and $\mu(F,\Psi) = \max_{1 \le i,j \le n} \left|(F\Psi)^i_j\right|$ is the mutual coherence between Ψ and F.*

The proof is mainly based on the concentration phenomenal of Lipschitz functions.

**Proof**

Let $x_0^T$ be the part of $x_0$ restricted on T, and

$$R(x_0) = \left\| \frac{1}{\sqrt{m}} U^\Omega x_0 \right\|^2 = \left\| \frac{1}{\sqrt{m}} U_T^\Omega x_0^T \right\|^2,$$

where $\|\ \|$ stands for the $\ell_2$-norm for a vector and the operator norm for a matrix. We expect that $R(x_0)$ concentrates around its expected value. For simplicity we suppose $\|x_0\| = \|x_0^T\| = 1$.

Note that since $E\{U_T^{\Omega*} U_T^\Omega\} = mI_T$

$$E\{R\} = \frac{1}{m} \text{Tr}(V^{\Omega*} V^\Omega) = \frac{1}{m} \sum_{k \in \Omega} \left\| x_0^{T*} \Psi_T^* \left(D^{(k-1)}\right)^* \right\|^2$$
$$= \frac{1}{m} \sum_{k \in \Omega} \|x_0^T\|^2 = 1$$

Assume $k \in \Omega$ and that $u_T^k = h^* D^{(k-1)} \Psi_T$ is the corresponding row of $U_T^\Omega$, so that

$$R = \frac{1}{m} \sum_{k \in \Omega} \left| \left\langle \Psi_T^* \left(D^{(k-1)}\right)^* h, x_0^T \right\rangle \right|^2$$

$$= \left\| \begin{bmatrix} x_0^{T*} \Psi_T^* \left(D^{(k_1-1)}\right)^* \\ x_0^{T*} \Psi_T^* \left(D^{(k_2-1)}\right)^* \\ \vdots \\ x_0^{T*} \Psi_T^* \left(D^{(k_m-1)}\right)^* \end{bmatrix} h \right\|^2 = \frac{1}{m} \|V^{\Omega*} h\|^2$$

where the notation $*$ means the conjugate transpose for a matrix such that

$$V^\Omega = \begin{bmatrix} D^{(k_1-1)} \Psi_T x_0^T & D^{(k_2-1)} \Psi_T x_0^T & \ldots & D^{(k_m-1)} \Psi_T x_0^T \end{bmatrix}$$

.

We now prove that when $x_0$ is fixed, $R^{1/2}$ is a Lipschitz function of the random vector **h**.

Let $f(\mathbf{h}) = R^{1/2} = \frac{1}{\sqrt{m}} \|V^{\Omega*} h\|$ and $\sigma_{V^\Omega} = \|V^\Omega\|$ for two independent copies of **h'**, **h"** of **h**,

$$|f(\mathbf{h'}) - f(\mathbf{h''})| = \frac{1}{\sqrt{m}} \left( \|V^\Omega \mathbf{h'}\| - \|V^\Omega \mathbf{h''}\| \right)$$

$$\le \frac{1}{\sqrt{m}} \|V^\Omega (\mathbf{h'} - \mathbf{h''})\| \le \frac{1}{\sqrt{m}} \sigma_{V^\Omega} \|\mathbf{h'} - \mathbf{h''}\|$$

.

where $f(\mathbf{h})$ is a $\frac{1}{\sqrt{m}} \sigma_{V^\Omega}$ - Lipschitz function of **h**.

As detailed in [18], [23], Lipschitz functions are very insensitive to local changes of the random vectors and are strongly concentrated around their means or medians. Respectively, we have the following tail bound for $f(h)$ when **h** is a Gaussian white vector,

$$\Pr\left(|f - m_f| > r\right) < \exp\left(-\frac{mr^2}{2\sigma_{V^\Omega}^2}\right), r > 0 \quad ,$$

where $m_f$ is the media of $f(\mathbf{h})$ used to product Gaussian measure. Since $m_f^2 = m_{f^2}$ and $|f - m_f| \le \sqrt{r} \Rightarrow |f^2 - m_f^2| \le r$ hold for $f, m_f \ge 0$, we get

$$\Pr\left(|R - m_R| > r\right) = \Pr\left(|f^2 - m_{f^2}| > r\right)$$

$$< \Pr\left(|f - m_f| > \sqrt{r}\right) < \exp\left(-\frac{mr}{2\sigma_{V^\Omega}^2}\right), r > 0$$

From [Theorem1.8, 20], we have

$$\Pr\left(|R - E\{R\}| > r + r_0\right) < \exp\left(-\frac{mr}{2\sigma_{V^\Omega}^2}\right), r > 0,$$

$$r_0 = \int_0^{+\infty} \exp\left(-\frac{mr}{2\sigma_{V^\Omega}^2}\right) dr = \frac{2\sigma_{V^\Omega}^2}{m}$$

.

so that

$$\Pr\left(|R - 1| > r + r_0\right) < \exp\left(-\frac{mr}{2\sigma_{V^\Omega}^2}\right), r > 0, r_0 = \frac{2\sigma_{V^\Omega}^2}{m}$$

That is,

$$\Pr(|R-1|>r) < e^2 \exp\left(-\frac{mr}{2\sigma_{V^\Omega}^2}\right), r > r_0 = \frac{2\sigma_{V^\Omega}^2}{m} \quad (3.2).$$

On the right hand side of the above inequality, $\sigma_{V^\Omega}$ is still a function of $x_0^T$. We now search an upper bound to $\sigma_{V^\Omega}$ by researching the property of the matrix $V^\Omega$. Since $V^\Omega$ is a submatrix of

$$V = \begin{bmatrix} l_1 & l_2 & & l_n \\ l_2 & & l_n & l_1 \\ & & & \\ l_n & l_1 & & l_{n-1} \end{bmatrix}$$

where $l = \Psi x_0^T = [l_1, l_2, ..., l_n]^T$ is a vector in the range of $\Psi_T$, we have

$$\sigma_{V^\Omega} = \|V^\Omega\| \le \|V\| = \sigma_V$$

Note that $V$ is similar to a convolution matrix. Let $F$ be the n-dim Discrete Fourier Transform matrix defined as

$$F = \left[F_k^j\right]_{n \times n}, F_k^j = e^{-i\frac{2\pi}{n}(j-1)(k-1)},$$

and so easily we have

$$V = F^* \text{diag}(\bar{l}) \text{conj}(F)/n,$$

where $\bar{l} = Fl$ is the Fourier transform of $l$.

As a result, $\sigma_V$ equals to the largest amplitude of $\bar{l}$, that is,

$$\sigma_V = \|Fl\|_\infty = \|F\Psi_T x_0^T\|_\infty$$

Let $\mu(F, \Psi) = \max_{1 \le i,j \le n} \left|(F\Psi)_j^i\right|$ be the coherence parameter between $F$ and $\Psi$ which is denoted as $\mu$ in short in the following paper.

We have

$$\sigma_V = \|F\Psi_T x_0^T\|_\infty = \left\|F \sum_{t \in T} x_0^t \psi_t\right\|_\infty$$
$$\le \sum_{t \in T} |x_0^t| \|F\psi_t\|_\infty \le \mu \sum_{t \in T} |x_0^t| \quad (3.3),$$
$$\le \mu \sqrt{\sum_{t \in T} |x_0^t|^2} \sqrt{S} = \mu \sqrt{S}$$

The first inequality follows from the triangle inequality while the last inequality follows from Cauchy's inequality. Taking (3.3) into (3.2) gives

$$\Pr\left(\left|\left\|\frac{1}{\sqrt{m}} U_T^\Omega x_0^T\right\|^2 - 1\right| > r\right) < e^2 \exp\left(-\frac{mr}{2S\mu^2}\right),$$
$$r > \frac{2S\mu^2}{m}$$

which establishes the claim. ∎

Lemma 3.1 does not assert that it holds for any signal with the same support. That's because when we arbitrarily choose a signal vector without concern for the specific sensing matrix, we just derive an a priori probable tail bound; however, there still exists some special signal vectors associated with the singular values of the sensing matrix which cannot be measured a priori. For example, let $g$ be an n-dimensional Gaussian white vector, and $f(g) = \max_{\|x\|=1} \langle g, x \rangle$. For arbitrarily chosen unit vector $x$ we have the following tail for $\langle g, x \rangle$ from [2.9, 20]

$$\Pr(|\langle g, x \rangle| > r) \le e^{-\frac{r^2}{2}}.$$

However, the above inequality does not hold for $f(g) = \|g\|$. The next theorem gives a more general concentration inequality than (3.1) and holds for any *signal $x_0$ support on T*, as well as for the RIP verification for our convolution matrix.

**Theorem 3.2**, *Let $\Psi$, H, U, $\Omega$ be the same as in Lemma 1and fix a subset T of size S on the signal domain, then for any signal $x_0$ support on T,*

$$\Pr\left(\left|\left\|\frac{1}{\sqrt{m}} U^\Omega x_0\right\|^2 - \|x_0\|^2\right| > r\|x_0\|^2\right)$$
$$< e^2 S^2 \exp\left(-\frac{mr}{4\mu^2 S}\right), for\ r > \frac{4\mu^2 S}{m} \quad (3.4)$$

*Indeed, if we unfix the support set T,*

$$\Pr\left(\left|\left\|\frac{1}{\sqrt{m}} U^\Omega x_0\right\|^2 - \|x_0\|^2\right| > r\|x_0\|^2\right)$$
$$< e^2 S^2 \binom{n}{S} \exp\left(-\frac{mr}{4\mu^2 S}\right), for\ r > \frac{4\mu^2 S}{m} \quad (3.5)$$

*holds for any signal $x_0$ with no more than S non-zero entries.*

The proof of this theorem is based on the following idea: when $R(x_0) = \left\|\frac{1}{\sqrt{m}} U^\Omega x_0\right\|^2 = \left\|\frac{1}{\sqrt{m}} U_T^\Omega x_0^T\right\|^2$ is bounded for a certain group of the S-dimensional unit vectors,

$R(x_0)$ is bounded on the whole S-dimensional unit ball.

**Proof**

Let $\|x_0\| = 1$, $x_0^T$ be the part of $x_0$ restricted on $T$, and consider $S$-dime vectors
$e_i = [0,...,0,1,0,...,0]^*, e_i(i)=1, e_i(\neq i)=0, 1 \leq i \leq n$,
$x_0^T$ can be expressed as $x_0^T = \sum_{i=1}^{S} c_i e_i$, with $\sum_{i=1}^{S} c_i^2 = 1$.

So
$$R(x_0) = \left\| \frac{1}{\sqrt{m}} U_T^\Omega x_0^T \right\|^2 = \left\| \frac{1}{\sqrt{m}} U_T^\Omega \sum_{i=1}^{S} c_i e_i \right\|^2$$
$$= \sum_{i=1}^{S} c_i^2 \left\| \frac{1}{\sqrt{m}} U_T^\Omega e_i \right\|^2 + \sum_{\substack{i,j=1 \\ i \neq j}}^{S} c_i c_j \left\langle \frac{1}{\sqrt{m}} U_T^\Omega e_i, \frac{1}{\sqrt{m}} U_T^\Omega e_j \right\rangle$$

Consider the vectors $(e_i - e_j)/\sqrt{2}, 1 \leq i,j \leq S, i \neq j$. Since
$$\left\langle \frac{1}{\sqrt{m}} U_T^\Omega e_i, \frac{1}{\sqrt{m}} U_T^\Omega e_j \right\rangle$$
$$= \frac{1}{2} \left( \left\| \frac{1}{\sqrt{m}} U_T^\Omega \frac{e_i - e_j}{\sqrt{2}} \right\|^2 - \frac{1}{2} \left\| \frac{1}{\sqrt{m}} U_T^\Omega e_i \right\|^2 - \frac{1}{2} \left\| \frac{1}{\sqrt{m}} U_T^\Omega e_i \right\|^2 \right)$$

If we have
$$1 - r \leq \left\| \frac{1}{\sqrt{m}} U_T^\Omega e_i \right\|^2, \left\| \frac{1}{\sqrt{m}} U_T^\Omega e_j \right\|^2,$$
$$\left\| \frac{1}{\sqrt{m}} U_T^\Omega \frac{e_i - e_j}{\sqrt{2}} \right\|^2 \leq 1 + r$$

we immediately get
$$-r \leq \left\langle \frac{1}{\sqrt{m}} U_T^\Omega e_i, \frac{1}{\sqrt{m}} U_T^\Omega e_j \right\rangle \leq r.$$

Then
$$R(x_0) = \sum_{i=1}^{S} c_i^2 \left\| \frac{1}{\sqrt{m}} U_T^\Omega e_i \right\|^2$$
$$+ \sum_{\substack{i,j=1 \\ i \neq j}}^{S} c_i c_j \left\langle \frac{1}{\sqrt{m}} U_T^\Omega e_i, \frac{1}{\sqrt{m}} U_T^\Omega e_j \right\rangle$$
$$\leq (1+r) \sum_{i=1}^{S} c_i^2 + \sum_{\substack{i,j=1 \\ i \neq j}}^{S} |c_i c_j| r$$
$$= 1 + \sum_{i,j=1}^{S} |c_i c_j| r \leq 1 + Sr$$

for $\sum_{i,j=1}^{S} |c_i c_j| = \left( \sum_{i=1}^{S} |c_i| \right)^2 \leq S$, and
$$R(x_0) \geq (1-r) \sum_{i=1}^{S} c_i^2 - \sum_{\substack{i,j=1 \\ i \neq j}}^{S} |c_i c_j| r$$
$$= 1 - \sum_{i,j=1}^{S} |c_i c_j| r \geq 1 - Sr$$

holds for any $x_0^T$ on the unit ball.

Now define some events
$$A_T = \left\{ \left| \left\| \frac{1}{\sqrt{m}} U_T^\Omega x_0^T \right\|^2 - 1 \right| \leq Sr \right\},$$
$$B_T^i = \left\{ \left| \left\| \frac{1}{\sqrt{m}} U_T^\Omega e_i \right\|^2 - 1 \right| \leq r \right\}, B_T = \bigcap_{1 \leq i \leq S} B_T^i,$$
$$C_T^{i,j} = \left\{ \left| \left\| \frac{1}{\sqrt{m}} U_T^\Omega \frac{e_i - e_j}{\sqrt{2}} \right\|^2 - 1 \right| \leq r \right\}, i \neq j, C_T = \bigcap_{\substack{i,j=1 \\ i \neq j}} C_T^{i,j}.$$

From the above discussion we know that when $B_T$ and $C_T$ happen, $A_T$ must happen. Accordingly,
$$\Pr(A_T^c) \leq \Pr(B_T^c \cup C_T^c) \leq \Pr(B_T^c) + \Pr(C_T^c)$$
$$\leq \sum_{1 \leq i \leq S} \Pr(B_T^{i\,c}) + \sum_{\substack{i,j=1 \\ i \neq j}}^{S} \Pr(C_T^{i,j\,c})$$

Since the supports of the vectors group $e_i, 1 \leq i \leq S$ and $(e_i - e_j)/\sqrt{2}, 1 \leq i,j \leq S, i \neq j$ are 1 and 2, respectively, from lemma 1 we have
$$\Pr(B_T^{i\,c}) = \Pr\left( \left| \left\| \frac{1}{\sqrt{m}} U_T^\Omega e_i \right\|^2 - 1 \right| > r \right)$$
$$< e^2 \exp\left( -\frac{mr}{2\mu^2} \right), r > \frac{2\mu^2}{m}$$
$$\Pr(C_T^{i,j\,c}, i \neq j) = \Pr\left( \left| \left\| \frac{1}{\sqrt{m}} U_T^\Omega \frac{e_i - e_j}{\sqrt{2}} \right\|^2 - 1 \right| > r \right)$$
$$< e^2 \exp\left( -\frac{mr}{4\mu^2} \right), r > \frac{4\mu^2}{m}$$

In conclusion

$$\Pr(A_T^c) = \Pr\left(\left|\left\|\frac{1}{\sqrt{m}}U_T^\Omega x_0^T\right\|^2 - 1\right| > Sr\right)$$

$$< e^2 S \exp\left(-\frac{mr}{2\mu^2}\right) + e^2 S(S-1)\exp\left(-\frac{mr}{4\mu^2}\right)$$

$$< e^2 S^2 \exp\left(-\frac{mr}{4\mu^2}\right), \text{ for } r > \frac{4\mu^2}{m}$$

Let $r \leftarrow Sr$, then

$$\Pr\left(\left|\left\|\frac{1}{\sqrt{m}}U_T^\Omega x_0^T\right\|^2 - 1\right| > r\right)$$

$$< e^2 S^2 \exp\left(-\frac{mr}{4\mu^2 S}\right), \text{ for } r > \frac{4\mu^2 S}{m}.$$

So for any signal $x_0$ with no more than S non-zero entries

$$\Pr\left(\left|\left\|\frac{1}{\sqrt{m}}U_T^\Omega x_0^T\right\|^2 - 1\right| > r\right)$$

$$< e^2 S^2 \binom{n}{S}\exp\left(-\frac{mr}{4\mu^2 S}\right), \text{ for } r > \frac{4\mu^2 S}{m}$$

The claim is then established.

∎

When (3.4) is established, the verification of RIP for a Gaussian white random convolution matrix is completed.

3.2 Proof of the main theorems for Gaussian ensembles

To establish the exact recovery theorem for our random convolution system followed by deterministic subsampling, we follow the program in [1]. As detailed in these references, the exact recovery of a signal $x_0$ supported on T with a given sign sequence $z$ from $y = U^\Omega x_0$ is performed by solving (1.2) if and only if there exists a *dual vector*

$$\pi = U^{\Omega*}U_T^\Omega\left(U_T^{\Omega*}U_T^\Omega\right)^{-1} z$$

such that $\pi(t) < 1$, $t \in T^C$, where $T^C$ is the complement of T in the signal domain. With (3.4) in theorem 3.2 established, we have

$$\Pr(\lambda_{\min} < 1-r \text{ or } \lambda_{\max} > 1+r)$$

$$< e^2 S^2 \exp\left(-\frac{mr}{4\mu^2 S}\right), \text{ for } r > \frac{4\mu^2 S}{m} \quad (3.6),$$

where $\lambda_{\max} = \sup_{\|x^T\|=1}\left\|\frac{1}{\sqrt{m}}U_T^\Omega x^T\right\|^2$ and

$\lambda_{\min} = \inf_{\|x^T\|=1}\left\|\frac{1}{\sqrt{m}}U_T^\Omega x^T\right\|^2$ are the largest and smallest eigenvalues of $\frac{1}{m}U_T^{\Omega*}U_T^\Omega$, *respectively.* That is,

$$\Pr\left(\left\|\frac{1}{m}U_T^{\Omega*}U_T^\Omega - I_T\right\| > r\right)$$

$$< e^2 S^2 \exp\left(-\frac{mr}{4\mu^2 S}\right), \text{ for } r > \frac{4\mu^2 S}{m} \quad (3.7)$$

We will prove theorem 1.1 with the help of the powerful inequalities (3.6)-(3.7).

Proof of Theorem 1.1:

Fix a $t_0 \in T^C$, suppose $w_{t_0} = \left(U_T^{\Omega*}U_T^\Omega\right)^{-1}U_T^{\Omega*}U_{t_0}^\Omega$, and $v_{t_0} = U_T^{\Omega*}U_{t_0}^\Omega$, then $\pi(t_0) = \langle w_{t_0}, z\rangle$. Following the program of [4], [6], we first derive the bound for the vector $v_{t_0}$.

As $v_{t_0} = U_T^{\Omega*}U_{t_0}^\Omega = \Psi_T^* H^{\Omega*}H^\Omega \psi_{t_0}$, and $\Psi_T^*\psi_{t_0} = \mathbf{0}$ for $\Psi$ is an orthobasis and $t_0 \notin T$,

$$v_{t_0} = U_T^{\Omega*}U_{t_0}^\Omega = \Psi_T^*\left(H^{\Omega*}H^\Omega - mI\right)\psi_{t_0}.$$

Let $T' = T \cup \{t_0\}$ and $v_{t_0}' = \Psi_{T'}^*\left(H^{\Omega*}H^\Omega - mI\right)\psi_{t_0}$.

Then $v_{t_0}$ is part of $v_{t_0}'$ by restricting $v_{t_0}'$ on T, and $v_{t_0}'$ is a column of

$$\Psi_{T'}^*\left(H^{\Omega*}H^\Omega - mI\right)\Psi_{T'} = m\left(\frac{1}{m}U_{T'}^{\Omega*}U_{T'}^\Omega - I_T\right).$$

So

$$\|v_{t_0}\| \leq \|v_{t_0}'\| \leq m\left\|\frac{1}{m}U_{T'}^{\Omega*}U_{T'}^\Omega - I_T\right\|$$

$$\|w_{t_0}\| = \left\|\left(U_{\Omega T}^*U_{\Omega T}\right)^{-1}v_{t_0}\right\| \leq \left\|\left(\frac{1}{m}U_{\Omega T}^*U_{\Omega T}\right)^{-1}\right\|\left\|\frac{1}{m}U_{T'}^{\Omega*}U_{T'}^\Omega - I\right\|$$

$$= \left\|\frac{1}{m}U_{T'}^{\Omega*}U_{T'}^\Omega - I\right\| \bigg/ \left\|\frac{1}{m}U_{\Omega T}^*U_{\Omega T}\right\|$$

Form (3.6)-(3.7) we get

$$\Pr\left(\begin{array}{l}1-r < \lambda_{\min}\left(\frac{1}{m}U_{T'}^{\Omega*}U_{T'}^\Omega\right) < \lambda_{\min}\left(\frac{1}{m}U_T^{\Omega*}U_T^\Omega\right)\\ < \lambda_{\max}\left(\frac{1}{m}U_T^{\Omega*}U_T^\Omega\right) < \lambda_{\max}\left(\frac{1}{m}U_{T'}^{\Omega*}U_{T'}^\Omega\right) < 1+r\end{array}\right)$$

$$> 1 - e^2(S+1)^2 \exp\left(-\frac{mr}{4\mu^2(S+1)}\right), \text{ for } r > \frac{4\mu^2(S+1)}{m}$$

for $|T'| = S+1$ and $T \subset T'$.

So

$$\Pr\left(\|w_{t_0}\| > \alpha\right) \le \Pr\left(\left\|\frac{1}{m}U^*_{\Omega T^{+1}}U_{\Omega T^{+1}} - I_{T^{+1}}\right\| > \alpha/(1+\alpha)\right)$$

$$< 1 - e^2(S+1)^2 \exp\left(-\frac{m\alpha/(1+\alpha)}{4\mu^2(S+1)}\right),$$

$$\text{for} \quad \alpha/(1+\alpha) > \frac{4\mu^2(S+1)}{m}$$

According to Hoeffeding's inequality [24], we get

$$\Pr\left(|\pi(t_0)| \ge 1\right) = \Pr\left(|\langle w_{t_0}, z\rangle| \ge 1\right) < 2\exp\left(-\frac{1}{2\|w_{t_0}\|^2}\right)$$

if $w_{t_0}$ is fixed. Accordingly,

$$\Pr\left(\sup_{t_0 \in T^C} |\pi(t_0)| \ge 1\right) < \Pr\left(\sup_{t_0 \in T^C} |\pi(t_0)| \ge 1 \Big| \sup_{t_0 \in T^C} \|w_{t_0}\| \le \alpha\right)$$

$$+ \Pr\left(\sup_{t_0 \in T^C} \|w_{t_0}\| > \alpha\right)$$

$$< n\Pr\left(|\pi(t_0)| \ge 1 \Big| \|w_{t_0}\| \le \alpha\right) + n\Pr\left(\|w_{t_0}\| > \alpha\right)$$

$$< 2n\exp\left(-\frac{1}{2\alpha^2}\right) + n\Pr\left(\|w_{t_0}\| > \alpha\right)$$

$$< 2n\exp\left(-\frac{1}{2\alpha^2}\right) + ne^2(S+1)^2 \exp\left(-\frac{m\alpha/(1+\alpha)}{4\mu^2(S+1)}\right),$$

$$\text{for} \quad \alpha/(1+\alpha) > \frac{4\mu^2(S+1)}{m}$$

We choose $\alpha = \sqrt{1/(2\log(4n/\delta))}$

such that $2n\exp\left(-\frac{1}{2\alpha^2}\right) = \frac{\delta}{2}$

For the second term

$$ne^2(S+1)^2 \exp\left(-\frac{m\alpha/(1+\alpha)}{4\mu^2(S+1)}\right) < \frac{\delta}{2},$$

we have

$$m > 4\mu^2(S+1)\left(1 + \sqrt{2\log(4n/\delta)}\right)$$
$$\cdot \log\left(2e^2(S+1)^2 n/\delta\right) \quad (3.8),$$

and for $\alpha/(1+\alpha) > \frac{4\mu^2(S+1)}{m}$, we have

$$m > 4\mu^2(S+1)\left(1 + \sqrt{2\log(4n/\delta)}\right)$$

which is weaker than (3.8) if $\delta$ is sufficiently small.

Choose a constant $K_1$ such that

$$1 + \sqrt{2\log(4n/\delta)} < K_1 \log(n/\delta)^{1/2},$$

then

$$m > 8K_1(S+1)\mu^2 \log(n/\delta)^{1/2}$$
$$\cdot \max\left(\log(2e^2(S+1)^2), \log(n/\delta)\right).$$

In conclusion, the exact recovery is ensured when the number of measurements $m$ obeys

$$m > K\mu^2 S \log(n/\delta)^{1/2}$$
$$\cdot \max\left(\log(2e^2(S+1)^2), \log(n/\delta)\right)$$

for a numerical constant K. The theory is proved. ∎

3.3 Symmetrical Bernoulli ensemble

When the rand waveform *h* is generated from a symmetrical Bernoulli distribution, that is,

$$\Pr(h_i = 1) = \Pr(h_i = -1) = \frac{1}{2}, 1 \le i \le n,$$

the proof of the main theorem still follows the program as detailed above. In this situation, the inequality (3.2) changes to

$$\Pr(|R-1| > r) < 2e^2 \exp\left(-\frac{mr}{16\sigma_{V^\Omega}^2}\right), r > r_0 = \frac{32\sigma_{V^\Omega}^2}{m}$$

(3.9).

Respectively, we have the following RIP verification theorem.

**Theorem 3.3** *Let $\Psi$, H, U, $\Omega$ be as in Lemma 1 except that h is generated from a symmetrical Bernoulli distribution, fix a subset T of size S on signal domain, then for any signal $x_0$ support on T,*

$$\Pr\left(\left|\left\|\frac{1}{\sqrt{m}}U^\Omega x_0\right\|^2 - \|x_0\|^2\right| > r\|x_0\|^2\right)$$
$$< 2e^2 S^2 \exp\left(-\frac{mr}{32\mu^2 S}\right), \text{for } r > \frac{64\mu^2 S}{m}$$

(3.10)

*Indeed, if unfix the support set T,*

$$\Pr\left(\left|\left\|\frac{1}{\sqrt{m}}U^\Omega x_0\right\|^2 - \|x_0\|^2\right| > r\|x_0\|^2\right)$$
$$< 2e^2 S^2 \binom{n}{S}\exp\left(-\frac{mr}{32\mu^2 S}\right), \text{for } r > \frac{64\mu^2 S}{m}$$

(3.11)

*holds for any signal $x_0$ with no more than S non-zero entries.*

We omit the details of proof for Theorem 3.3 as well as the theorem 1.2.

## IV. Conclusion

In this paper we analyze the CS convolution framework, convolving the tested signal with a white random waveform, followed by subsampling at fixed locations in the measurement domain, i.e., equal interval sampling. Discussions are limited to circular convolution where linear convolution can be easily transformed to circular convolution. As an effect of the reduction in the freedom of randomness, the linear CS convolution system needs more measurements than the circular one with the same size. It also becomes inefficient when the waveform length is too short relative to the signal length. Indeed, in some applications the bandwidth of the random waveform is shorter than the tested signal. However, even in this case one may have super-resolution results when the original signal is sparse enough or when more prior information is used. Such super-resolution effects are beyond the scope of this paper and are the subjects of current research.

ACKNOWLEDGEMENT

## References

[1] E. Candès, J. Romberg, and T. Tao, "Robust uncertainty principles: Exact signal reconstruction from highly incomplete frequency information," *IEEE Trans. on Information Theory*, 52(2) pp. 489 - 509, February 2006

[2] E. Candès and T. Tao, "Near optimal signal recovery from random projections: Universal encoding strategies?", *IEEE Trans. on Information Theory*, 52(12), pp. 5406 - 5425, December 2006

[3] E. Candès, J. Romberg, and T. Tao, Stable signal recovery from incomplete and inaccurate measurements. Communications on Pure and Applied Mathematics, 59(8), pp. 1207-1223, August 2006

[4] E. Candès and J. Romberg, Sparsity and incoherence in compressive sampling. (Inverse Problems, 23(3) pp. 969-985, 2007)

[5] W. Bajwa, J. Haupt, G. Raz, S. Wright, and R. Nowak, Toeplitz-structured compressed sensing matrices. (IEEE Workshop on Statistical Signal Processing (SSP), Madison, Wisconsin, August 2007)

[6] J. Romberg, Compressive sensing by random convolution. (Preprint, 2008) Available: http://dsp.rice.edu/files/cs/RandomConvolution.pdf

[7] M. Duarte, M. Davenport, D. Takhar, J. Laska, T. Sun, K. Kelly, and R. Baraniuk, Single-pixel imaging via compressive sampling. (IEEE Signal Processing Magazine, 25(2), pp. 83 - 91, March 2008)

[8] L. Jacques, P. Vandergheynst, A. Bibet, V. Majidzadeh, A. Schmid, and Y. Leblebici, CMOS compressed imaging by random comvolution. (Preprint, 2008)

[9] J. Tropp and A. Gilbert, Signal recovery from random measurements via orthogonal matching pursuit. (IEEE Trans. on Information Theory, 53(12) pp. 4655-4666, December 2007).

[10] J. Tropp, M. Wakin, M. Duarte, D. Baron, and R. Baraniuk, Random filters for compressive sampling and reconstruction. (IEEE Int. Conf. on Acoustics, Speech, and Signal Processing (ICASSP), Toulouse, France, May 2006)

[11] S. R. J. Axelsson, Random noise radar/sodar with ultrawideband waveforms, IEEE Trans. Geosci. Remote Sens., 45 (2007), pp. 1099-1114.

[12] W. U. Bajwa, J. D. Haupt, G. M. Raz, S. J. Wright, and R. D. Nowak, Toeplitz-structured compressed sensing matrices, in Proc. IEEE Stat. Sig. Proc. Workshop, Madison, WI, August 2007, pp. 294-298.

[13] R. Baraniuk and P. Steeghs, Compressive radar imaging, in Proc. IEEE Radar Conference, Boston, MA, April 2007, pp. 128-133.

[14] R. G. Baraniuk, M. Davenport, R. DeVore, and M. Wakin, A simple proof of the restricted isometry property for random matrices. to appear in Constructive Approximation, 2008.

[15] D. L. Donoho, Compressed sensing, IEEE Trans. Inform. Theory, 52 (2006), pp. 1289-1306.

[16] D.L. Donoho, Yaakov Tsaig, Iddo Drori, and Jean-Luc Starck, Sparse solution of underdetermined linear equations by stagewise orthogonal matching pursuit. (Preprint, 2007)

[17] M. A. Herman and T. Strohmer, High-resolution radar via compressed sensing. Submitted to IEEE. Trans. Sig. Proc., 2008.

[18] M. Ledoux, The Concentration of Measure Phenomenon, American Mathematical Society, 2001.

[19] M. Richards, Fundamentals of Radar Signal Processing, McGraw-Hill, 2005.

[20] R. Marcia, Z. Harmany, R. Willett, Compressive Coded Aperture Imaging. (SPIE Electronic Imaging, 2009).

[21] R. Baraniuk and P. Steeghs, Compressive radar imaging. (IEEE Radar Conference, Waltham, Massachusetts, April 2007)

[22] S. Bhattacharya, T. Blumensath, B. Mulgrew, and M. Davies, Fast encoding of synthetic aperture radar raw data using compressed sensing. (IEEE Workshop on Statistical Signal Processing, Madison, Wisconsin, August 2007)

[23] M. Talagrand, A New Look at independence, Annal. Prob., vol.24, no. 1, 1-34(1996).

[24] W. Hoeffding, Probability inequalities for sums of bounded random variables, Journal of the American Statistical Association 58 (301): 13–30, March 1963.

[25] Radu Berinde and Piotr Indyk, Sparse recovery using sparse random matrices. (Preprint, 2008)
Available: http://people.csail.mit.edu/indyk/report.pdf

[26] Sami Kirolos, Jason Laska, Michael Wakin, Marco Duarte, Dror Baron, Tamer Ragheb, Yehia Massoud, and Richard Baraniuk, Analog-to-information conversion via random demodulation. (IEEE Dallas Circuits and Systems Workshop (DCAS), Dallas, Texas, 2006)

[27] D. Achlioptas, Database-friendly random projections, Proc. ACM SIGACT-SIGMOD-SIGART Symp. on Principles of Database Systems (2001), pp. 274–281.

[28] H. Rauhut, Circulant and Toeplitz matrices in compressed sensing, Proc. SPARS 2009.